\documentclass[12pt]{article}
\usepackage{amsmath,amsfonts,amssymb}

\newtheorem{thm}{Theorem}[section]

\newtheorem{lem}[thm]{Lemma}

\def\qed{\nopagebreak\hfill{\rule{4pt}{7pt}}}
\def\pf{\noindent {\it Proof.} }
\def\ps{\mathrm{ps}}
\begin{document}

\begin{center}
{\large  The $q$-Log-convexity of the Generating Functions of the
Squares of Binomial Coefficients}
\end{center}

\begin{center}
William Y. C. Chen$^{1}$, Robert L. Tang$^2$, \\
Larry X. W. Wang$^{3}$ and
Arthur L. B. Yang$^{4}$\\[6pt]
Center for Combinatorics, LPMC-TJKLC\\
Nankai University, Tianjin 300071, P. R. China\\[5pt]
$^{1}${\tt chen@nankai.edu.cn}, $^{2}${\tt
tangling@cfc.nankai.edu.cn}, $^{3}${\tt wxw@cfc.nankai.edu.cn},
$^{4}${\tt yang@nankai.edu.cn}
\end{center}

\vspace{0.3cm} \noindent{\bf Abstract.} We prove a conjecture of Liu
and Wang on the $q$-log-convexity of the polynomial sequence
$\{\sum_{k=0}^n{n\choose k}^2q^k\}_{n\geq 0}$. By using Pieri's rule
and the Jacobi-Trudi identity for Schur functions, we obtain an
expansion of a sum of products of elementary symmetric functions in
terms of Schur functions with nonnegative coefficients. Then the
principal specialization leads to the $q$-log-convexity. We also
prove that a technical condition of Liu and Wang  holds for the
squares of the binomial coefficients. Hence we deduce that the
linear transformation with respect to the triangular array
$\{{n\choose k}^2\}_{0\leq k\leq n}$ is log-convexity preserving.



%%%\noindent {\bf Suggested Running Title:} $q$-log-convexity

\allowdisplaybreaks

\section{Introduction}

The  objective of this paper is to prove a conjecture of Liu and
Wang \cite{wang-liu} on the $q$-log-convexity of the following
polynomials
\begin{equation}
W_n(q)=\sum_{k=0}^n{n\choose k}^2q^k.
\end{equation}
The polynomial $W_n(q)$ has appeared as the rank generating function
of the lattice of noncrossing partitions of type $B$ on
 $[n]$, see Reiner \cite{reiner1997}. For the type $A$ case,
the rank generating  function of the lattice of noncrossing
partitions on $[n]$ is equal to the Narayana polynomial
\begin{equation}
N_n(q)=\sum_{k=0}^n \frac{1}{1}{n \choose k}{n\choose k+1}q^k.
\end{equation}
Liu and Wang \cite{wang-liu} also conjectured that $N_n(q)$ are
$q$-log-convex. Subsequently, this conjecture was proved by Chen,
Wang and Yang \cite{Chen-Wang-Yang}.

The polynomials $W_n(q)$ also arise in the theory of growth series
of the root lattice. Recall that the classical root lattice $A_n$ is
generated by $\mathcal{M}=\{{\bf e}_i-{\bf e}_j: 0\leq i,j\leq n+1
\, {\rm with} \, i\neq j\}$. Then the growth series is defined to be
the generating function
$$G(q)=\sum_{k\geq 0}S(k)q^k,$$
where $S(k)$ is the number of elements ${\bf u}\in A_n$ with length
$k$. It is known that $G(q)$ is a rational function of the form
$$G(q)=\frac{h(q)}{(1-q)^d},$$
where $d$ is the rank of $A_n$ and $h(q)$ is a polynomial of degree
less than or equal to $d$. The polynomial $h(q)$ is defined to be
the coordinator polynomial of the growth series, see \cite{Ben}.
Recently, Ardila et al. \cite{Ardila} have shown that the above
coordinator polynomial $h(q)$ of $A_n$ equals the polynomial
$W_n(q)$.

Recall that a sequence $\{a_k\}_{k\geq 0}$ of nonnegative numbers is
log-convex if $a_{k-1}a_{k+1}\geq a_{k}^2$ for any $k\geq 1$. The
$q$-log-convexity is a property defined in \cite{wang-liu} for
sequences of polynomials over the field of real numbers, similar to
the concept of $q$-log-concavity introduced by Stanley and
subsequently studied by Butler \cite{butler1990}, Krattenthaler
\cite{kratte1989}, Leroux \cite{leroux1990} and Sagan
\cite{sagan1992}. Given  a sequence $\{f_n(q)\}_{n\geq 0}$, we say
that it is \emph{$q$-log-convex} if for any $k\geq 1$ the difference
$$f_{k+1}(q)f_{k-1}(q)-f_k(q)^2$$ has nonnegative coefficients as a
polynomial of $q$. It has been shown that many combinatorial
polynomials are $q$-log-convex, such as the Bell polynomials, the
Eulerian polynomials, the Bessel polynomials, the Ramanujan
polynomials and the Dowling polynomials, see Liu and Wang
\cite{wang-liu}, and Chen, Wang and Yang \cite{Chen-Wang-Yang1}.

Clearly, if the sequence $\{f_n(q)\}_{n\geq 0}$ is $q$-log-convex,
then for each fixed positive number $q$ the sequence
$\{f_n(q)\}_{n\geq 0}$ is log-convex. It can be easily verified that
the sequence of the central binomial coefficients $\{b_n\}_{n\geq
0}$ is log-convex, where
$$b_n={2n \choose n} =\sum_{k=0}^n{n\choose k}^2.$$
Do\v{s}li\'{c} and Veljan \cite{doslic} obtained the log-convexity
of the sequence of the central Delannoy numbers $\{d_n\}_{n\geq 0}$,
where
$$d_n=\sum_{k=0}^n{n\choose k}^22^k.$$ These two examples lead to the
conjecture of the $q$-log-convexity of $W_n(q)$, see
\cite[Conjecture 5.3]{wang-liu}. The first result of this paper is
to give an affirmative answer to this conjecture.

\begin{thm}\label{Liuwang1}
The polynomials $W_n(q)$ form a $q$-log-convex sequence.
\end{thm}

 By using the principal
specialization,  the $q$-log-convexity of $W_n(q)$ follows from the
Schur positivity of a sum of products of elementary symmetric
functions. This Schur positivity is based on an identity on
symmetric functions, which is proved by showing that both sides
satisfy the same recurrence relations. To establish the required
recurrence relation, we employ the Jacobi-Trudi identity and Pieri's
rule for Schur functions. We would like to note that in general the
polynomials
$$\sum_{k=0}^n{n\choose
k}^mq^k, \quad n\geq 0$$ are not $q$-log-convex for $m\geq 3$.

The second result of this paper is concerned with the condition on
linear transformations that preserve log-convexity.  Let
$\{a(n,k)\}_{0\leq k\leq n}$ be a triangular array of numbers. The
linear transformation on a sequence $\{x_k\}_{k\geq 0}$ with respect
to a triangular array $\{a(n,k)\}_{0\leq k\leq n}$ is defined by
\[ y_n=\sum_{k=0}^na(n,k)x_k.\]
Such a transformation is called \emph{log-convexity preserving} if
$\{y_n\}_{n\geq 0}$ is log-convex whenever $\{x_k\}_{k\geq 0}$ is
log-convex. Liu and Wang \cite{wang-liu} obtained a sufficient
condition under which the transformation with respect to a given
triangular array is log-convexity preserving.

Given a triangular array $\{a(n,k)\}_{0\leq k\leq n}$, define
$\alpha(n,r,k)$ by
\begin{equation*}
\alpha(n,r,k)=a(n+1,k)a(n-1,r-k)+a(n+1,r-k)a(n-1,k)
-2a(n,r-k)a(n,k),
\end{equation*}
where $n\geq 1$, $0\leq r\leq 2n$ and $0\leq k\leq \lfloor
\frac{r}{2}\rfloor$. The sufficient condition of Liu and Wang is
stated as follows.

\begin{thm}[{\cite[Theorem 4.8]{wang-liu}}]
\label{Liuwang3}  Assume that the polynomials
\[A_n(q)=\sum_{k=0}^na(n,k)q^k\] form a $q$-log-convex sequence. For
any given $n$ and $r$, if there exists an integer $k'=k'(n,r)$ such
that $\alpha(n,r,k)\geq 0$ for $k\leq k'$ and $\alpha(n,r,k)\leq 0$
for $k> k'$, then the linear transformation with respect to the
triangular array $\{a(n,k)\}_{0\leq k\leq n}$ is log-convexity
preserving.
\end{thm}

Liu and Wang conjectured that the above sufficient condition holds
for the coefficients of $W_n(q)$. We will also prove this
conjecture. So we deduce the following conclusion.

\begin{thm}\label{Liuwang2}
The linear transformation with respect to the triangular array
$\{{n\choose k}^2\}_{0\leq k\leq n}$ is log-convexity preserving.
\end{thm}

This paper is organized as follows. We  recall some definitions and
known results on symmetric functions in Section 2. In Section 3, we
will give an induction proof of the identity which implies the
desired Schur positivity, though a little more complicated. In
Section 4, we will give  the proofs of Theorem \ref{Liuwang1} and
Theorem \ref{Liuwang2}.

\section{Background on symmetric functions}

Throughout this paper we will adopt notation and terminology on
partitions and symmetric functions in Stanley \cite{stanley1999}.
Recall that a \emph{partition} $\lambda$ of a nonnegative integer
$n$ is a weakly decreasing sequence $(\lambda_1, \lambda_2, \ldots)$
of nonnegative integers satisfying $\sum_{i}\lambda_i=n$, denoted
$\lambda\vdash n$. We usually omit the parts $\lambda_i=0$. We also
denote a partition $\lambda\vdash n$ by
$(n^{m_n},\ldots,2^{m_2},1^{m_1})$ if $\lambda$ has $m_i$ $i$'s for
$1\leq i\leq n$. Let $\mathrm{Par}(n)$ denote the set of all
partitions of $n$.

If $\lambda\vdash n$, we draw a left-justified array of $n$ squares
with $\lambda_i$ squares in the $i$-th row. This array is called the
\emph{Young diagram} of $\lambda$. By transposing the diagram of
$\lambda$, we get the \emph{conjugate partition} of $\lambda$,
denoted $\lambda'$. We use  $\mu\subseteq \lambda$ to denote that
the Young diagram of $\mu$ is contained in the diagram of $\lambda$.

A semistandard Young tableau of shape $\lambda$ is an array
$T=(T_{ij})$ of positive integers of shape $\lambda$ such that it is
weakly increasing in each row and strictly increasing in each
column. The type of $T$ is defined as the composition
$\alpha=(\alpha_1,\,\alpha_2,\ldots)$, where $\alpha_i$ is the
number of $i$'s in $T$. Let $x$ denote the variables $\{
x_1,x_2,\ldots\}$. If $\mathrm{type}(T)=\alpha$, then we write
$$x^T=x_1^{\alpha_1}x_2^{\alpha_2}\cdots.$$
The \emph{Schur function} $s_{\lambda}(x)$ is defined as the
generating function
$$s_{\lambda}(x)=\sum_T x^T,$$
summed over all semistandard Young tableaux $T$ of shape $\lambda$.
When $\lambda=\emptyset$, we set $s_{\emptyset}(x)=1$.

It is well known that  Schur functions $s_{\lambda}(x)$ form a basis
for the ring of symmetric functions. A symmetric function $f(x)$ is
called \emph{Schur positive} if the coefficients $a_\lambda$ are all
nonnegative in the Schur expansion
$f(x)=\sum_{\lambda}a_{\lambda}s_{\lambda}(x)$.

When $\lambda=(1^k)$ for $k\geq 1$, the Schur function
$s_{\lambda}(x)$ becomes the $k$-th elementary symmetric function
$e_k(x)$, i.e.,
\begin{equation}\label{ele-def}
s_{(1^k)}(x)=e_k(x)=\sum_{1\leq i_1<\cdots<i_k}x_{i_1}\cdots
x_{i_k}.
\end{equation}

The dual Jacobi-Trudi identity gives an expression of the Schur
function $s_{\lambda}(x)$ in terms of elementary symmetric
functions.

\begin{thm}[{\cite[Corollary 7.16.2]{stanley1999}}]\label{jacobi}
Let $\lambda$ be a partition with the largest part $\leq n$ and
$\lambda'$ its conjugate. Then
$$s_\lambda(x)=\det(e_{{\lambda_i}'-i+j}(x))_{i,j=1}^n,$$
where $e_0=1$ and $e_k=0$ for $k<0$.
\end{thm}

Given any symmetric function $f(x)$, we may omit the variable set
$x$ if no confusion arises in the context. Now let us review the
definition of the principal specialization $\ps_n^1$ of a symmetric
function. For any symmetric function $f$, the action of $\ps_n^1$ is
defined as
$$\ps_n^1(f)=f(\underbrace{1,\ldots,1}_{n\, 1's},0,0,\ldots).$$
In particular, by \eqref{ele-def}, we have $\ps_n^1(e_k)={n\choose
k}$ and
\begin{equation}\label{decOrder}
\ps_n^1(e_k)  =  \ps_{n-1}^1(e_k+e_{k-1}),
\end{equation}
which is a restatement of the relation
\begin{equation*}
\binom{n}{k}  = \binom{n-1}{k}+\binom{n-1}{k-1}.
\end{equation*}

We will also  need the dual version of Pieri's rule which  expresses
the product of a Schur function $s_{\mu}$ and an elementary
symmetric function $e_k$ in terms of Schur functions.

\begin{thm}[{\cite{stanley1999}}]\label{pieri}
We have
$$s_{\mu}e_k=\sum_{\lambda}s_{\lambda},$$
summed over all partitions $\lambda$ such that $\mu\subseteq
\lambda$ and the difference of the Young diagrams of $\lambda$ and
$\mu$ contains no two squares in the same row.
\end{thm}

\section{A Schur Positivity Identity}

The objective of this section is to establish the following Schur
positivity theorem which will  be used to prove the
$q$-log-convexity of $W_n(q)$.

\begin{thm}\label{condc1}
For any $r\geq 1$, we have
\begin{equation}\label{main-eq}
{\sum_{k=0}^r
(e_{k-1}e_{k-1}e_{r-k}e_{r-k}+e_{k-2}e_{k}e_{r-k}e_{r-k}-2e_{k-1}e_{k}e_{r-k-1}e_{r-k})}=
\sum_{\lambda}s_{\lambda},
\end{equation}
where $\lambda$ sums over all partitions of $2r-2$ of the form
$(4^{i_4},3^{2i_3},2^{2i_2},1^{2i_1})$ with $i_1,i_2,i_3,i_4$ being
nonnegative integers.
\end{thm}

Before proving the above theorem, it is informative to give examples
for $r=3,4,5$. By  using the Maple package  ACE \cite{veigneau1998},
or SF \cite{stembridge1995}, we find that
\begin{align*}
&{\sum_{k=0}^3
(e_{k-1}e_{k-1}e_{3-k}e_{3-k}+e_{k-2}e_{k}e_{3-k}e_{3-k}-2e_{k-1}e_{k}e_{3-k-1}e_{3-k})}\\
&\quad=s_{(1^4)}+s_{(2^2)}+s_{(4)},\\
&{\sum_{k=0}^4
(e_{k-1}e_{k-1}e_{4-k}e_{4-k}+e_{k-2}e_{k}e_{4-k}e_{4-k}-2e_{k-1}e_{k}e_{4-k-1}e_{4-k})}\\
&\quad=s_{(1^6)}+s_{(2^2,1^2)}+s_{(4,1^2)}+s_{(3^2)}, \\
&{\sum_{k=0}^5
(e_{k-1}e_{k-1}e_{5-k}e_{5-k}+e_{k-2}e_{k}e_{5-k}e_{5-k}-2e_{k-1}e_{k}e_{5-k-1}e_{5-k})}\\
&\quad=s_{(4,2^2)}+s_{(4^2)}+s_{(1^8)}+s_{(2^2,1^4)}+s_{(2^4)}+s_{(4,1^4)}+s_{(3^2,1^2)}.
\end{align*}

Let $L(r)$ and $R(r)$ denote the left-hand side and the right-hand
side of \eqref{main-eq}, respectively. The key idea of the proof of
Theorem \ref{condc1} is to show that $L(r)$ and $R(r)$ satisfy the
same recurrence relations. It is easy to find the recurrence
relation of $R(r)$. To derive the recurrence relation of $L(r)$, we
will give several lemmas. For the sake of presentation, we will use
the following notation.  For $t\geq 0$, let
\begin{align*}
A_1(t,k,i,j)&=e_{k} s_{(3^i,2^{k-i+j},1^{4t-3k-2j-i})},
\\[3pt]
A_2(t,k,i,j)&=e_{k}s_{(3^i,2^{k-i+j},1^{4t-3k-2j-i})},
\\[3pt]
A_3(t,k,i,j)&=e_{k-1}s_{(3^i,2^{k-i+j},1^{4t-3k-2j-i+1})},\\[3pt]
A_4(t,k,i,j)& =e_{k}s_{(3^i,2^{k-i+j-1},1^{4t-3k-2j-i+2})},\\[3pt]
B_1(t,k,i,j)&=e_{k} s_{(3^i,2^{k-i+j},1^{4t-3k-2j-i-2})},\\[3pt]
B_2(t,k,i,j)&=e_{k} s_{(3^i,2^{k-i+j},1^{4t-3k-2j-i-2})},\\[3pt]
B_3(t,k,i,j)&=e_{k-1} s_{(3^i,2^{k-i+j},1^{4t-3k-2j-i-1})},\\[3pt]
B_4(t,k,i,j)&=e_{k} s_{(3^i,2^{k-i+j-1},1^{4t-3k-2j-i})},
\end{align*}
and  let
\begin{align*}
A_1(t)&=\sum_{k=0}^t \sum_{i=0}^{k}\sum_{j=0}^{2t-2k}A_1(t,k,i,j),
\\[3pt]
A_2(t)&=\sum_{k=0}^{t-1}
\sum_{i=0}^{k}\sum_{j=0}^{2t-2k-1}A_2(t,k,i,j),
\\[3pt]
A_3(t)&= \sum_{k=1}^{t} \sum_{i=0}^{k}\sum_{j=0}^{2t-2k}A_3(t,k,i,j),\\[3pt]
A_4(t)& =  \sum_{k=1}^t \sum_{i=0}^{k-1}\sum_{j=0}^{2t-2k+1}A_4(t,k,i,j),\\[3pt]
B_1(t)& = \sum_{k=0}^{t-1} \sum_{i=0}^{k}\sum_{j=0}^{2t-2k-1}
B_1(t,k,i,j),\\[3pt]
B_2(t)&= \sum_{k=0}^{t-1} \sum_{i=0}^{k}\sum_{j=0}^{2t-2k-2}
B_2(t,k,i,j),\\[3pt]
B_3(t)&=\sum_{k=1}^{t-1} \sum_{i=0}^{k}\sum_{j=0}^{2t-2k-1}
B_3(t,k,i,j),\\[3pt]
B_4(t)&=\sum_{k=1}^t \sum_{i=0}^{k-1}\sum_{j=0}^{2t-2k}B_4(t,k,i,j).
\end{align*}

The following lemma gives explicit expressions for $L(r)$ according
to the parity of $r$.

\begin{lem}\label{left}
For $t\geq 0$ we have
\begin{align*}
L(2t+1) & =  A_1(t)+A_2(t)-A_3(t)-A_4(t), \\
L(2t) & =  B_1(t)+B_2(t)-B_3(t)-B_4(t).
\end{align*}
\end{lem}

\pf For any $r\geq 0$, we have
\begin{align*}
L(r)  = & \sum_{k=0}^r
(e_{k-1}e_{k-1}e_{r-k}e_{r-k}+e_{k-2}e_{k}e_{r-k}e_{r-k}-2e_{k-1}e_{k}e_{r-k-1}e_{r-k})\\
  = & \sum_{k=0}^r
(e_{k-1}e_{k-1}e_{r-k}e_{r-k}-e_{k-1}e_{k}e_{r-k-1}e_{r-k})\\
   & \quad+ \sum_{k=0}^r
(e_{k-2}e_{k}e_{r-k}e_{r-k}-e_{k-1}e_{k}e_{r-k-1}e_{r-k})\\
 = & \sum_{k=0}^r e_{k-1}e_{r-k} \det\begin{pmatrix}e_{k-1} & e_{k}\\
e_{r-k-1} & e_{r-k}\end{pmatrix} \\
& + \, \sum_{k=0}^r e_{k}e_{r-k} \det\begin{pmatrix}e_{k-2} & e_{k-1}\\
e_{r-k-1} & e_{r-k}
\end{pmatrix}.
\end{align*}
By the dual Jacobi-Trudi identity, we obtain that
\begin{align*}
L(r)  = & \sum_{{k=0 \atop k-1\geq r-k}}^r e_{k-1}e_{r-k}s_{(2^{r-k},1^{2k-r-1})} - \sum_{{k=0 \atop k-1< r-k-1}}^r e_{k-1}e_{r-k}s_{(2^k,1^{r-2k-1})}\\[5pt]
& + \, \sum_{{k=0 \atop k-2\geq r-k}}^r
e_{k}e_{r-k}s_{(2^{r-k},1^{2k-r-2})} - \sum_{{k=0 \atop k-2<
r-k-1}}^r e_{k}e_{r-k}s_{(2^{k-1},1^{r-2k})}.
\end{align*}
Applying the dual version of Pieri's rule to the products
$e_{k-1}s_{(2^{r-k},1^{2k-r-1})}$, $e_{r-k}s_{(2^k,1^{r-2k-1})},
e_{k-1}s_{(2^{r-k},1^{2k-r-2})}$ and
$e_{r-k}s_{(2^{k-1},1^{r-2k})}$, we get
\begin{align*}
L(r)  = & \sum_{{k=0 \atop k-1\geq r-k}}^r
\sum_{i=0}^{r-k}\sum_{j=0}^{2k-r-1} e_{r-k}
s_{(3^i,2^{r-k-i+j},1^{2k-r-2-j+k-i-j})}\\
&- \, \sum_{{k=0 \atop k-1< r-k-1}}^r \sum_{i=0}^{k}\sum_{j=0}^{r-2k-1}e_{k-1}s_{(3^i,2^{k-i+j},1^{r-2k-j+r-k-i-j-1})}\\
& + \, \sum_{{k=0 \atop k-2\geq r-k}}^r
\sum_{i=0}^{r-k}\sum_{j=0}^{2k-r-2} e_{r-k}
s_{(3^i,2^{r-k-i+j},1^{2k-r-2-j+k-i-j})}\\
& -\, \sum_{{k=0 \atop k-2< r-k-1}}^r
\sum_{i=0}^{k-1}\sum_{j=0}^{r-2k} e_{k}
s_{(3^i,2^{k-i+j-1},1^{r-2k-j+r-k-i-j})}.
\end{align*}
Setting $r=2t$ or $r=2t+1$, we get the required relations. \qed

To find a recurrence relation of $L(r)$, let us recall an operator
$\Delta^{\mu}$ associated with a partition $\mu$,  which acts on
symmetric functions. This operator was introduced by Chen, Wang and
Yang \cite{Chen-Wang-Yang}. Given two partitions $\lambda$ and
$\mu$, let $\lambda\cup\mu$ be the partition whose parts are
obtained by taking the union of the parts of $\lambda$ and $\mu$.
For a symmetric function $f$ with the expansion
$$f=\sum_{\lambda}a_{\lambda}s_{\lambda},$$ the action of
$\Delta^{\mu}$ on $f$ is defined by
$$\Delta^{\mu}(f)=\sum_{{\lambda}}a_{\lambda}s_{\lambda\cup\mu}.$$

In order to compute the difference $L(2t+1)-\Delta^{(1,1)}L(2t)$, we
need to evaluate $A_m(t,k,i,j) - \Delta^{(1,1)}(B_m(t,k,i,j))$ for
$1\leq m\leq 4$. In fact, we will be able to express these
differences as double sums of Schur functions. For $t\geq 0$, let
\begin{align*}
T_1(t,i,j,k) &=\sum\limits_{a=0}^{\beta_1+1} \quad
\sum\limits_{b=0}^{\min(\beta_2,\beta_3+1)}P(t,k+1,i,j-1,a,b),\\[5pt]
T_2(t,i,j,k) &=\sum\limits_{a=0}^{\beta_1} \quad
\sum\limits_{b=0}^{\min(\beta_2,\beta_3)}P(t,k,i,j,a,b),\\[5pt]
T_3(t,i,j,k) &=\sum\limits_{a=0}^{\beta_1+1} \quad
\sum\limits_{b=0}^{\min(\beta_2,\beta_3+1)}P(t,k+1,i,j-1,a,b),\\[5pt]
T_4(t,i,j,k) &=\sum\limits_{a=0}^{\beta_1} \quad
\sum\limits_{b=0}^{\min(\beta_2,\beta_3)}P(t,k,i,j,a,b),\\[5pt]
T_5(t,i,j,k) &=\sum\limits_{a=0}^{\beta_1-1}
\sum\limits_{b=0}^{\min(\beta_2,\beta_3-1)}P(t,k,i,j,a,b),\\[5pt]
T_6(t,i,j,k)&=\sum\limits_{a=0}^{\beta_1-2} \quad
\sum\limits_{b=0}^{\min(\beta_2,\beta_3-2)}P(t,k-1,i,j+1,a,b),\\[5pt]
T_7(t,i,j,k)& =\sum\limits_{a=0}^{\beta_1-1} \quad
\sum\limits_{b=0}^{\min(\beta_2-1,\beta_3-1)}P(t,k,i,j,a,b), \\[5pt]
T_8(t,i,j,k) &=\sum\limits_{a=0}^{\beta_1-2} \quad
\sum\limits_{b=0}^{\min(\beta_2-1,\beta_3-2)}P(t,k-1,i,j+1,a,b),
\end{align*}
where
$$
\begin{array}{rcl}
P(t,k,i,j,a,b)&=&s_{(4^a,3^{i-a+b},2^{4t-2k-2i-j-b},1^{4k+i+2j-a-b-4t})},\\[5pt]
\beta_1 &=& 4k+i+2j-4t,\\[5pt]
\beta_2 &=& k-i+j, \\[5pt]
\beta_3 &=& 4k+i+2j-a-4t.
\end{array}
$$

Then we have the following result.

\begin{lem}\label{ab2t}
Suppose that $t\geq 0$.
\begin{itemize}
\item[({i})] If $0\leq k\leq t-1, 0\leq i\leq k, 0\leq j\leq 2t-2k-1$,
then
$$
A_1(t,k,i,j) -
\Delta^{(1,1)}(B_1(t,k,i,j))=T_1(t,i,j,k)+T_2(t,i,j,k).
$$

\item[({ii})] If $0\leq k\leq t-1, 0\leq i\leq k, 0\leq j\leq 2t-2k-2$, then
$$
A_2(t,k,i,j) -
\Delta^{(1,1)}(B_2(t,k,i,j))=T_3(t,i,j,k)+T_4(t,i,j,k).
$$

\item[({iii})] If $1\leq k\leq t-1, 0\leq i\leq k, 0\leq j\leq 2t-2k-1$, then
$$
A_3(t,k,i,j) -
\Delta^{(1,1)}(B_3(t,k,i,j))=T_5(t,i,j,k)+T_6(t,i,j,k).
$$

\item[({iv})] If $1\leq k\leq t, 0\leq i\leq k-1, 0\leq j\leq 2t-2k$, then
$$
A_4(t,k,i,j) -
\Delta^{(1,1)}(B_4(t,k,i,j))=T_7(t,i,j,k)+T_8(t,i,j,k).
$$
\end{itemize}
\end{lem}

\pf We only give the proof of (i). The proofs of ({ii}),({iii}) and
({iv}) are  analogous to that of ({i}) and are omitted. Recall that
$$A_1(t,k,i,j)=e_{k} s_{(3^i,2^{k-i+j},1^{4t-3k-2j-i})},\, B_1(t,k,i,j)=e_{k} s_{(3^i,2^{k-i+j},1^{4t-3k-2j-i-2})}.$$
Using the dual version of Pieri's rule, we obtain that
$$A_1(t,k,i,j)=\sum_{(a,b,c,d)}s_{(4^a,3^{i-a+b},2^{k-i+j-b+c},1^{4t-3k-2j-i-c+d})},$$
summed over all nonnegative integer sequences $(a,b,c,d)$ satisfying
$$a+b+c+d=k,\, a \leq i,\, b \leq k-i+j,\, c\leq 4t-3k-2j-i.$$
Note that the shape
$(4^a,3^{i-a+b},2^{k-i+j-b+c},1^{4t-3k-2j-i-c+d})$ is obtained from
the Young diagram of $(3^i,2^{k-i+j},1^{4t-3k-2j-i})$ by adding $a$
squares in the fourth column, $b$ squares in the third column, $c$
squares in the second column, and $d$ squares in the first column.
Similarly, we see that
$$B_1(t,k,i,j)=\sum_{(a,b,c,d)}s_{(4^a,3^{i-a+b},2^{k-i+j-b+c},1^{4t-3k-2j-i-2-c+d})},$$
summed over all nonnegative integer sequences $(a,b,c,d)$ satisfying
$$a+b+c+d=k,\,a \leq i, \,b \leq k-i+j,\, c\leq 4t-3k-2j-i-2.$$
Therefore,
\begin{align*}
&A_1(t,k,i,j) - \Delta^{(1,1)}(B_1(t,k,i,j))\\[5pt]
&=\sum\limits_{{(a,b,c,d)\atop
c=4t-3k-2j-i}}s_{(4^a,3^{i-a+b},2^{k-i+j-b+c},1^{4t-3k-2j-i-c+d})}\\[5pt]
 &\qquad+\sum\limits_{{(a,b,c,d)\atop
c=4t-3k-2j-i-1}}s_{(4^a,3^{i-a+b},2^{k-i+j-b+c},1^{4t-3k-2j-i-c+d})},
\end{align*}
where both sums range over nonnegative integers $a,b,d$ satisfying
$$a \leq i,\, b \leq k-i+j,\, a+b+c+d=k.$$
Since $0\leq j\leq 2t-2k-1$, we have
\begin{eqnarray*}
\lefteqn{k-(4t-3k-2j-i)=4k+2j+i-4t}&\\
&\rule{100pt}{0pt}\leq 4k+2(2t-2k-1)+i-4t\leq i,
\end{eqnarray*}
and
\begin{eqnarray*}
\lefteqn{k-(4t-3k-2j-i-1)=4k+2j+i-4t+1}&\\
&\rule{100pt}{0pt}\leq 4k+2(2t-2k-1)+i-4t+1\leq i.
\end{eqnarray*}
Thus
$$ A_1(t,k,i,j) -
\Delta^{(1,1)}(B_1(t,k,i,j))=T_1(i,j,k)+T_2(i,j,k).
$$
This completes the proof of ({i}). \qed

In light of the above lemma, we will show that
$L(2t+1)-\Delta^{(1,1)}L(2t)$ can be expressed in terms of ten parts
$T_1(t), T_2(t), \ldots, T_{10}(t)$, as defined below,
\begin{align*}
T_1(t)= & \sum_{k=0}^{t-1} \sum_{i=0}^{k}\sum_{j=0}^{2t-2k-1}T_1(t,i,j,k),\quad T_2(t)=  \sum_{k=0}^{t-1} \sum_{i=0}^{k}\sum_{j=0}^{2t-2k-1}T_2(t,i,j,k),\\
T_3(t)= & \sum_{k=0}^{t-1}
\sum_{i=0}^{k}\sum_{j=0}^{2t-2k-2}T_3(t,i,j,k), \quad T_4(t)=  \sum_{k=0}^{t-1} \sum_{i=0}^{k}\sum_{j=0}^{2t-2k-2}T_4(t,i,j,k),\\
T_5(t)= & \sum_{k=1}^{t-1}
\sum_{i=0}^{k}\sum_{j=0}^{2t-2k-1}T_5(t,i,j,k),\quad T_6(t)=  \sum_{k=1}^{t-1} \sum_{i=0}^{k}\sum_{j=0}^{2t-2k-1}T_6(t,i,j,k),\\
T_7(t)= &
\sum_{k=1}^{t}\sum_{i=0}^{k-1}\sum_{j=0}^{2t-2k}T_7(t,i,j,k),\quad
\rule{10pt}{0pt} T_8(t)=
\sum_{k=1}^{t}\sum_{i=0}^{k-1}\sum_{j=0}^{2t-2k}T_8(t,i,j,k), \\
T_9(t)= & \sum_{k=0}^t e_ks_{(3^k,2^{2t-2k})}, \quad
\rule{50pt}{0pt} T_{10}(t)=
\sum_{k=0}^{t-1}e_ks_{(3^{k+1},2^{2t-2k-2},1)}.
\end{align*}

\begin{lem}\label{lem-ten} For $t\geq 0$, we have
\begin{eqnarray*}
\lefteqn{L(2t+1)-\Delta^{(1,1)}L(2t)=(T_1(t)+T_2(t)+T_3(t)+T_4(t))}&\\[3pt]
&\rule{100pt}{0pt}-(T_5(t)+T_6(t)+T_7(t)+T_8(t))+(T_9(t)-T_{10}(t)).
\end{eqnarray*}
\end{lem}
\pf By Lemma \ref{left},  $L(2t+1)-\Delta^{(1,1)}L(2t)$ equals
\begin{align*}
 (A_1(t)&-\Delta^{(1,1)}B_1(t))+(A_2(t)-\Delta^{(1,1)}B_2(t))\\[3pt]
             & - \,
             (A_3(t)-\Delta^{(1,1)}B_3(t))-(A_4(t)-\Delta^{(1,1)}B_4(t))\\[3pt]
           = & \left(\sum_{k=0}^t \sum_{i=0}^{k}\sum_{j=0}^{2t-2k}A_1(t,k,i,j) -
           \sum_{k=0}^{t-1}
\sum_{i=0}^{k}\sum_{j=0}^{2t-2k-1}
\Delta^{(1,1)}(B_1(t,k,i,j))\right)\\[6pt]
& + \, \left(\sum_{k=0}^{t-1}
\sum_{i=0}^{k}\sum_{j=0}^{2t-2k-1}A_2(t,k,i,j) - \sum_{k=0}^{t-1}
\sum_{i=0}^{k}\sum_{j=0}^{2t-2k-2} \Delta^{(1,1)}(B_2(t,k,i,j))\right)\\[6pt]
& - \, \left(\sum_{k=1}^{t}
\sum_{i=0}^{k}\sum_{j=0}^{2t-2k}A_3(t,k,i,j) - \sum_{k=1}^{t-1}
\sum_{i=0}^{k}\sum_{j=0}^{2t-2k-1} \Delta^{(1,1)}(B_3(t,k,i,j))\right)\\[6pt]
& - \, \left(\sum_{k=1}^t
\sum_{i=0}^{k-1}\sum_{j=0}^{2t-2k+1}A_4(t,k,i,j) - \sum_{k=1}^t
\sum_{i=0}^{k-1}\sum_{j=0}^{2t-2k}
\Delta^{(1,1)}(B_4(t,k,i,j))\right)\\
           = & \sum_{k=0}^t \sum_{i=0}^{k} A_1(t,k,i,2t-2k)+\sum_{k=0}^{t-1} \sum_{i=0}^{k} A_2(t,k,i,2t-2k-1)\\
             & + \, \sum_{k=0}^{t-1} \sum_{i=0}^{k}\sum_{j=0}^{2t-2k-1}\left(A_1(t,k,i,j) - \Delta^{(1,1)}(B_1(t,k,i,j))\right)\\
             & + \, \sum_{k=0}^{t-1} \sum_{i=0}^{k}\sum_{j=0}^{2t-2k-2}\left(A_2(t,k,i,j) - \Delta^{(1,1)}(B_2(t,k,i,j))\right)\\
             & - \, \sum_{k=1}^t \sum_{i=0}^{k} A_3(t,k,i,2t-2k)-\sum_{k=1}^{t} \sum_{i=0}^{k-1} A_4(t,k,i,2t-2k+1)\\
             & - \, \sum_{k=1}^{t-1} \sum_{i=0}^{k}\sum_{j=0}^{2t-2k-1}\left(A_3(t,k,i,j) - \Delta^{(1,1)}(B_3(t,k,i,j))\right)\\
             & - \, \sum_{k=1}^{t} \sum_{i=0}^{k-1}\sum_{j=0}^{2t-2k}\left(A_4(t,k,i,j) -
             \Delta^{(1,1)}(B_4(t,k,i,j))\right).
\end{align*}
Applying Lemma \ref{ab2t}, we get the desired relation. \qed

Although Lemma \ref{lem-ten} gives an expression of
$L(2t+1)-\Delta^{(1,1)}L(2t)$, it is difficult to derive a
recurrence relation of $L(r)$ since each $T_i(t)$ is a quintuple sum
for $1\leq i \leq 8$. We still need to compute
\begin{equation*}
(L(2(t+1)+1)-\Delta^{(1,1)}L(2(t+1)))-\Delta^{(2,2)}(L(2t+1)-\Delta^{(1,1)}L(2t)).
\end{equation*}
By Lemma \ref{lem-ten}, the above expression equals
\begin{align}
&\sum_{i=1}^4(T_i(t+1)-\Delta^{(2,2)}T_i(t))-\sum_{i=5}^8(T_i(t+1)-\Delta^{(2,2)}T_i(t))\nonumber\\
&\rule{20pt}{0pt}+(T_{9}(t+1)-\Delta^{(2,2)}T_{9}(t))
+(T_{10}(t+1)-\Delta^{(2,2)}T_{10}(t)).\label{2-recrr}
\end{align}
For notational convenience, let
\begin{align}
Q_1(t,k,i,j,a,b) &=
s_{(4^a,3^{k-a+b+1},2^{4t-4k-j-b-1},1^{5k+2j-a-b-4t+3})},\label{nota-q}\\[5pt]
Q_2(t,k,i,j,a,b) &=
s_{(4^a,3^{k-a+j+1},2^{4t-3k-i-2j},1^{3k+2i+j-a-4t+1})},\\[5pt]
\gamma_1&= 5k+2j+2-4t,\\[5pt]
\gamma_2&= 5k+2j-a+2-4t,\label{nota-2}\\[5pt]
\gamma_3&= 3k+2i+j-4t,
\end{align}

and let
\begin{align*}
T_{11}(t)=&T_{21}(t)=\sum_{k=0}^{t-1}\sum_{j=0}^{2t-2k-1}\sum_{a=0}^{\gamma_1}\sum_{b=0}^{\min(j,\gamma_2)} Q_1(t,k,i,j,a,b),\\[5pt]
T_{12}(t)=&T_{22}(t)=\sum\limits_{k=0}^t\sum\limits_{i=0}^k\quad
\sum\limits_{j=0}^{2t-2k-1}\quad \sum\limits_{a=0}^{\gamma_3}
Q_2(t,k,i,j,a,b),\\[5pt]
T_{31}(t)=&\sum_{k=0}^t\sum_{j=0}^{2t-2k-2}\sum_{a=0}^{\gamma_1}\sum_{b=0}^{\min(j,\gamma_2)} Q_1(t,k,i,j,a,b),\\[5pt]
T_{32}(t)=&\sum\limits_{k=0}^t\sum\limits_{i=0}^k\quad
\sum\limits_{j=0}^{2t-2k-2}\quad \sum\limits_{a=0}^{\gamma_3}
Q_2(t,k,i,j,a,b),\\[5pt]
T_{41}(t)=&\sum_{k=0}^{t-1}\sum_{j=-1}^{2t-2k-3}\sum_{a=0}^{\gamma_1+1}\sum_{b=0}^{\min(j+1,\gamma_2+1)} Q_1(t,k,i,j,a,b),\\[5pt]
T_{42}(t)=&\sum\limits_{k=0}^{t-1}\sum\limits_{i=-1}^{k-1}\quad
\sum\limits_{j=0}^{2t-2k-2}\quad \sum\limits_{a=0}^{\gamma_3+1}
Q_2(t,k,i,j,a,b),\\[5pt]
T_{51}(t)=&\sum_{k=0}^{t-1}\sum_{j=-1}^{2t-2k-2}\sum_{a=0}^{\gamma_1}\sum_{b=0}^{\min(j+1,\gamma_2)} Q_1(t,k,i,j,a,b),\\[5pt]
T_{52}(t)=&\sum\limits_{k=0}^{t-1}\sum\limits_{i=-1}^{k-1}\quad
\sum\limits_{j=0}^{2t-2k-1}\quad \sum\limits_{a=0}^{\gamma_3}
Q_2(t,k,i,j,a,b),\\[5pt]
T_{61}(t)=&\sum_{k=0}^{t-1}\sum_{j=-2}^{2t-2k-3}\sum_{a=0}^{\gamma_1+1}\sum_{b=0}^{\min(j+2,\gamma_2+1)} Q_1(t,k,i,j,a,b),\\[5pt]
T_{62}(t)=&\sum\limits_{k=-2}^{t-3}\sum\limits_{i=0}^{k+2}\quad
\sum\limits_{j=2}^{2t-2k-3}\quad \sum\limits_{a=0}^{\gamma_3+1}
Q_2(t,k,i,j,a,b),\\[5pt]
T_{71}(t)=&\sum_{k=0}^{t}\sum_{j=-2}^{2t-2k-2}\sum_{a=0}^{\gamma_1+1}\sum_{b=-1}^{\min(j+1,\gamma_2)} Q_1(t,k,i,j,a,b),\\[5pt]
T_{72}(t)=&\sum\limits_{k=0}^{t}\sum\limits_{i=0}^{k-1}\quad
\sum\limits_{j=-1}^{2t-2k-1}\quad \sum\limits_{a=0}^{\gamma_3}
Q_2(t,k,i,j,a,b),\\[5pt]
T_{81}(t)=&\sum_{k=0}^{t-1}\sum_{j=0}^{2t-2k-2}\sum_{a=0}^{\gamma_1+1}\sum_{b=0}^{\min(j,\gamma_2+1)} Q_1(t,k,i,j,a,b),\\[5pt]
T_{82}(t)=&\sum\limits_{k=-1}^{t-1}\sum\limits_{i=0}^{k}\quad
\sum\limits_{j=0}^{2t-2k-2}\quad \sum\limits_{a=0}^{\gamma_3+1}
Q_2(t,k,i,j,a,b).
\end{align*}

The following lemma shows that each $T_i(t+1)-\Delta^{(2,2)}T_i(t)$
can be divided into two parts for $1\leq i\leq 8$, with each part
being a quadruple sum of Schur functions.

\begin{lem}\label{T1-T8}
For $t\geq 0$ and $1\leq i\leq 8$, we have
\begin{equation*}
T_i(t+1)-\Delta^{(2,2)}T_i(t)= T_{i1}(t)+T_{i2}(t).
\end{equation*}
\end{lem}

\pf We will give  only the proof of  the identity for $T_1(t)$,
since the other cases can be verified by the same argument. Observe
that
\begin{align*}
T_1(t)=\sum_{k=0}^{t-1} \sum_{i=0}^{k}\sum_{j=0}^{2t-2k-1}\quad
\sum\limits_{a=0}^{4k+i+2j-4t+1} \quad
\sum\limits_{b=0}^{\min(k-i+j,4k+i+2j-a-4t+1)}\\[6pt]
\left(s_{(4^a,3^{i-a+b},2^{4t-2k-2i-j-b-1},1^{4k+i+2j-a-b-4t+2})}\right).\rule{30pt}{0pt}
\end{align*}
It follows that
\begin{align*}
T_1(t+1)&=\sum_{k=0}^{t} \sum_{i=0}^{k}\sum_{j=0}^{2t-2k+1}\quad
\sum\limits_{a=0}^{4k+i+2j-4t-3} \quad
\sum\limits_{b=0}^{\min(k-i+j,4k+i+2j-a-4t-3)}\\[6pt]
&\rule{50pt}{0pt}\left(s_{(4^a,3^{i-a+b},2^{4t-2k-2i-j-b+3},1^{4k+i+2j-a-b-4t-2})}\right)\\[10pt]
&=\sum_{k=-1}^{t-1} \sum_{i=0}^{k+1}\sum_{j=0}^{2t-2k-1}\quad
\sum\limits_{a=0}^{4k+i+2j-4t+1} \quad
\sum\limits_{b=0}^{\min(k-i+j+1,4k+i+2j-a-4t+1)}\\[6pt]
&\rule{50pt}{0pt}\left(s_{(4^a,3^{i-a+b},2^{4t-2k-2i-j-b+1},1^{4k+i+2j-a-b-4t+2})}\right)\\[10pt]
&=\sum_{k=0}^{t-1} \sum_{i=0}^{k+1}\sum_{j=0}^{2t-2k-1}\quad
\sum\limits_{a=0}^{4k+i+2j-4t+1} \quad
\sum\limits_{b=0}^{\min(k-i+j+1,4k+i+2j-a-4t+1)}\\[6pt]
&\rule{50pt}{0pt}\left(s_{(4^a,3^{i-a+b},2^{4t-2k-2i-j-b+1},1^{4k+i+2j-a-b-4t+2})}\right),
\end{align*}
where the last equality holds because the upper bound $4k+i+2j-4t+1$
of $a$ is negative for $k=-1$. Hence we deduce that
\begin{align*}
&T_1(t+1)-\Delta^{(2,2)}T_1(t)\\[6pt]
&=\sum_{k=0}^{t-1} \sum_{i=k+1}^{k+1}\sum_{j=0}^{2t-2k-1}\quad
\sum\limits_{a=0}^{4k+i+2j-4t+1} \quad
\sum\limits_{b=0}^{\min(k-i+j+1,4k+i+2j-a-4t+1)}\\[6pt]
&\rule{50pt}{0pt}\left(s_{(4^a,3^{i-a+b},2^{4t-2k-2i-j-b+1},1^{4k+i+2j-a-b-4t+2})}\right)\\[10pt]
&\rule{10pt}{0pt} + \sum_{k=0}^{t-1}
\sum_{i=0}^{k}\sum_{j=0}^{2t-2k-1}\quad
\sum\limits_{a=0}^{4k+i+2j-4t+1} \quad
\sum\limits_{b=0}^{\min(k-i+j+1,4k+i+2j-a-4t+1)}\\[6pt]
&\rule{50pt}{0pt}\left(s_{(4^a,3^{i-a+b},2^{4t-2k-2i-j-b+1},1^{4k+i+2j-a-b-4t+2})}\right)\\[10pt]
&\rule{10pt}{0pt} - \sum_{k=0}^{t-1}
\sum_{i=0}^{k}\sum_{j=0}^{2t-2k-1}\quad
\sum\limits_{a=0}^{4k+i+2j-4t+1} \quad
\sum\limits_{b=0}^{\min(k-i+j,4k+i+2j-a-4t+1)}\\[6pt]
&\rule{50pt}{0pt}\left(s_{(4^a,3^{i-a+b},2^{4t-2k-2i-j-b+1},1^{4k+i+2j-a-b-4t+2})}\right)\\[10pt]
&=\sum_{k=0}^{t-1} \sum_{i=k+1}^{k+1}\sum_{j=0}^{2t-2k-1}\quad
\sum\limits_{a=0}^{4k+i+2j-4t+1} \quad
\sum\limits_{b=0}^{\min(k-i+j+1,4k+i+2j-a-4t+1)}\\[6pt]
&\rule{50pt}{0pt}\left(s_{(4^a,3^{i-a+b},2^{4t-2k-2i-j-b+1},1^{4k+i+2j-a-b-4t+2})}\right)\\[10pt]
&\rule{10pt}{0pt} + \sum_{k=0}^{t-1}
\sum_{i=0}^{k}\sum_{j=0}^{2t-2k-1}\quad \sum\limits_{{a=0\atop
k-i+j+1\leq 4k+i+2j-a-4t+1 }}^{4k+i+2j-4t+1} \quad
\sum\limits_{b=0}^{k-i+j+1}\\[6pt]
&\rule{50pt}{0pt}\left(s_{(4^a,3^{i-a+b},2^{4t-2k-2i-j-b+1},1^{4k+i+2j-a-b-4t+2})}\right)\\[10pt]
&\rule{10pt}{0pt} - \sum_{k=0}^{t-1}
\sum_{i=0}^{k}\sum_{j=0}^{2t-2k-1}\quad \sum\limits_{{a=0\atop
k-i+j+1\leq 4k+i+2j-a-4t+1 }}^{4k+i+2j-4t+1} \quad
\sum\limits_{b=0}^{k-i+j}\\[6pt]
&\rule{50pt}{0pt}\left(s_{(4^a,3^{i-a+b},2^{4t-2k-2i-j-b+1},1^{4k+i+2j-a-b-4t+2})}\right)\\[10pt]
&=\sum_{k=0}^{t-1} \sum_{i=k+1}^{k+1}\sum_{j=0}^{2t-2k-1}\quad
\sum\limits_{a=0}^{4k+i+2j-4t+1} \quad
\sum\limits_{b=0}^{\min(k-i+j+1,4k+i+2j-a-4t+1)}\\[6pt]
&\rule{50pt}{0pt}\left(s_{(4^a,3^{i-a+b},2^{4t-2k-2i-j-b+1},1^{4k+i+2j-a-b-4t+2})}\right)\\[10pt]
&\rule{10pt}{0pt} + \sum_{k=0}^{t-1}
\sum_{i=0}^{k}\sum_{j=0}^{2t-2k-1}\quad \sum\limits_{{a=0\atop
k-i+j+1\leq 4k+i+2j-a-4t+1 }}^{4k+i+2j-4t+1} \quad
\sum\limits_{b=k-i+j+1}^{k-i+j+1}\\[6pt]
&\rule{50pt}{0pt}\left(s_{(4^a,3^{i-a+b},2^{4t-2k-2i-j-b+1},1^{4k+i+2j-a-b-4t+2})}\right)\\[10pt]
&=\sum_{k=0}^{t-1}\sum_{j=0}^{2t-2k-1}\sum_{a=0}^{5k+2j+2-4t}\sum_{b=0}^{\min(j,5k+2j-a+2-4t)}\\[6pt]
&\rule{50pt}{0pt} \left(s_{(4^a,3^{k-a+b+1},2^{4t-4k-j-b-1},1^{5k+2j-a-b-4t+3})}\right)\\[5pt]
&\rule{10pt}{0pt} +\sum\limits_{k=0}^{t-1}\sum\limits_{i=0}^k\quad
\sum\limits_{j=0}^{2t-2k-1}\quad \sum\limits_{a=0}^{3k+2i+j-4t}\\[6pt]
&\rule{50pt}{0pt}
\left(s_{(4^a,3^{k-a+j+1},2^{4t-3k-i-2j},1^{3k+2i+j-a-4t+1})}\right),
\end{align*}
as desired. This completes the proof. \qed

To compute \eqref{2-recrr}, it is necessary to simplify
$$\sum_{i=1}^4(T_i(t+1)-\Delta^{(2,2)}T_i(t))-\sum_{i=5}^8(T_i(t+1)-\Delta^{(2,2)}T_i(t)),$$
by the above lemma, which equals
$$\sum_{i=1}^4 (T_{i1}(t)+T_{i2}(t))-\sum_{i=5}^8 (T_{i1}(t)+T_{i2}(t)).$$
Moreover, we need to group the terms to reduce the relevant
quadruple sums to triple sums, and then to double sums. For $t\geq
0$, let
\begin{align*}
N_1(t)=&\sum_{k=0}^{t-1}\sum_{j=0}^{2t-2k-1}\sum_{a=0}^{5k+j-4t-1}
s_{(4^a,3^{k-a+j+1},2^{4t-4k-2j+1},1^{5k+j-1-4t-a})},\\[5pt]
N_2(t)=&\sum_{k=0}^{t-1}\sum_{j=0}^{2t-2k-2}\sum_{a=0}^{5k+j-4t+2}
s_{(4^a,3^{k-a+j+1},2^{4t-4k-2j-1},1^{5k+j+3-4t-a})},\\[5pt]
N_3(t)=&\sum_{k=0}^{t-1}\sum_{a=0}^{k}\sum_{b=0}^{\min(2t-2k-2,k-a)}s_{(4^a,3^{k-a+b+1},2^{2t-2k-b},1^{k-a-b+1})},\\[5pt]
N_4(t)=&\sum_{k=0}^{t-1}\sum_{a=0}^{k}\sum_{b=0}^{\min(2t-2k-1,k-a)}
s_{(4^a,3^{k+1-a+b},2^{2t-2k-b},1^{k+1-a-b})},\\[5pt]
N_5(t)=&\sum_{k=0}^{t-1}\sum_{j=0}^{2t-2k-2}\sum_{a=\max(0,5k+j-4t+3)}^{5k+2j-4t+3}
s_{(4^a,3^{6k-4t+2j-2a+4},2^{8t-9k-3j+a-4})},\\[5pt]
N_6(t)=&\sum_{k=0}^{t-1}\sum_{j=0}^{2t-2k-1}\sum_{a=\max(0,5k+j-4t+1)}^{5k+2j-4t+1}
s_{(4^a,3^{6k-4t+2j-2a+2},2^{8t-9k-3j+a-1})},\\[5pt]
M_1(t)=&\sum\limits_{k=0}^{t-1}\quad
\sum\limits_{j=0}^{2t-2k-2}\quad \sum\limits_{a=0}^{5k+j+2-4t}
s_{(4^a,3^{k+1-a+j},2^{4t-2j-1-4k},1^{5k+j-4t+3-a})},\\[5pt]
M_2(t)=&\sum\limits_{k=0}^{t-1}\quad
\sum\limits_{j=0}^{2t-2k-2}\quad \sum\limits_{a=0}^{5k+j-4t-1}
s_{(4^a,3^{k+1-a+j},2^{4t+1-4k-2j},1^{5k+j-4t-1-a})},\\[5pt]
M_3(t)=&\sum\limits_{k=0}^{t-1}\sum\limits_{i=1}^{k}\quad
 \sum\limits_{a=0}^{3k+j+2i-4t-3}s_{(4^a,3^{2t-k-a},2^{k+4-i},1^{k+2i-2t-4-a})},\\[5pt]
M_4(t)=&\sum\limits_{k=0}^{t-1}\sum\limits_{i=0}^k\quad
\sum\limits_{a=0}^{k+2i-2t-1}
s_{(4^a,3^{2t-a-k},2^{k-i+2},1^{k+2i-2t-a})},\\[5pt]
M_5(t)=&\sum\limits_{k=0}^{t-1}\sum\limits_{i=0}^k\quad
\sum\limits_{j=0}^{2t-2k-2}\sum\limits_{a=\max(0,3k+j+2i-4t+1)}^{3k+j+2i-4t+1}
s_{(4^{a},3^{k+1-a+j},2^{4t-2j-i-3k})},\\[5pt]
M_6(t)=&\sum_{k=1}^{t-1}\sum_{i=0}^k\sum_{j=0}^{2t-2k-1}\sum_{a=\max(0,3k+j+2i-4t-1)}^{3k+j+2i-4t-1}
s_{(4^{a},3^{k+1-a+j},2^{4t-2j-i-3k+1})}.
\end{align*}

The following lemma gives a strategy to group the terms of
$T_{i1}(t)$ and $T_{i2}(t)$, which leads to the reduction from
quadruple sums to triple sums.

\begin{lem}\label{t2nm} For $t\geq 0$, we have
\begin{align*}
T_{41}(t)-T_{61}(t)=&-N_1(t),\\[5pt]
T_{31}(t)-T_{71}(t)=&N_2(t)-N_3(t),\\[5pt]
T_{11}(t)-T_{81}(t)=&N_4(t)-N_5(t),\\[5pt]
T_{21}(t)-T_{51}(t)=&N_6(t),\\[5pt]
T_{32}(t)-T_{72}(t)=& -M_1(t)\\[5pt]
T_{42}(t)-T_{62}(t)=&M_2(t)-M_3(t),\\[5pt]
T_{12}(t)-T_{82}(t)=&M_4(t)-M_5(t),\\[5pt]
T_{22}(t)-T_{52}(t)=&M_6(t).
\end{align*}
\end{lem}
\pf We will prove the first identity, since the others can be proved
analogously. Using the notation  $Q_1(t,k,i,j,a,b)$, $ \gamma_1$ and
$\gamma_2$, as given in (\ref{nota-q})-(\ref{nota-2}),  we find
\begin{align*}
T_{41}(t)-T_{61}(t)=& \sum_{k=0}^{t-1}\sum_{j=-1}^{2t-2k-3}\sum_{a=0}^{\gamma_1+1}\sum_{b=0}^{\min(j+1,\gamma_2+1)} Q_1(t,k,i,j,a,b)\\[5pt]
              & -\sum_{k=0}^{t-1}\sum_{j=-2}^{2t-2k-3}\sum_{a=0}^{\gamma_1+1}\sum_{b=0}^{\min(j+2,\gamma_2+1)}
              Q_1(t,k,i,j,a,b)\\[5pt]
              =& \sum_{k=0}^{t-1}\sum_{j=-2}^{2t-2k-3}\sum_{a=0}^{\gamma_1+1}\sum_{b=0}^{\min(j+1,\gamma_2+1)} Q_1(t,k,i,j,a,b)\\[5pt]
              & -\sum_{k=0}^{t-1}\sum_{j=-2}^{2t-2k-3}\sum_{a=0}^{\gamma_1+1}\sum_{b=0}^{\min(j+2,\gamma_2+1)}
              Q_1(t,k,i,j,a,b)\\[5pt]
              = & -\sum_{k=0}^{t-1}\sum_{j=-2}^{2t-2k-3}\sum_{{a=0 \atop j+2\leq \gamma_2+1}}^{\gamma_1+1}\sum_{b=j+2}^{j+2}
              Q_1(t,k,i,j,a,b)\\[5pt]
              = & -\sum_{k=0}^{t-1}\sum_{j=-2}^{2t-2k-3}\sum_{{a=0}}^{5k+j-4t+1}
              Q_1(k,i,j,a,j+2).
\end{align*}
Substituting $j$ by $j-2$ in the last summation, we are led to the
required relation. This completes the proof. \qed

For $t\geq 0$, let
\begin{align*}
C_1(t)=&\sum_{k=1}^{t}\sum_{a=0}^{5k-4t-3}s_{(4^a,3^{k-a},2^{4t-4k+3},1^{5k-4t-a-2})},\\[5pt]
C_2(t)=&\sum_{k=1}^{t-1}\sum_{j=0}^{2t-2k-1}\sum_{a=\max(0,5k+j-4t-1)}^{5k+j-4t-1}s_{(4^a,3^{k-a+j+1},2^{4t-4k-2j+1})},\\[5pt]
C_3(t)=&\sum_{k=0}^{t-1}\sum_{a=0}^{3k-2t+1}s_{(4^a,3^{2t-k-a},2,1^{3k-2t-a+2})},\\[5pt]
C_4(t)=&\sum_{k=0}^{t-1}\sum_{a=\max(0,5k-4t+1)}^{5k-4t+1}s_{(4^a,3^{k+1-a},2^{4t-4k})},\\[5pt]
C_5(t)=&\sum_{k=1}^{t-1}\sum_{j=0}^{2t-2k-2}\sum_{a=\max(0,5k+j-4t+2)}^{5k+j-4t+2}s_{(4^{a},3^{k+j-a+2},2^{4t-4k-2j-2})},\\[5pt]
D_1(t)=&\sum_{k=0}^{t-2}\sum_{j=1}^{2t-2k-3}\sum_{a=\max(0,5k+j-4t+3)}^{5k+j-4t+3}s_{(4^a,3^{k+1-a+j},2^{4t-2j-4k-1})},\\[5pt]
D_2(t)=&\sum_{a=0}^{t-3}s_{(4^a,3^{t-a},2^3,1^{t-a-2})},\\[5pt]
D_3(t)=&\sum_{k=0}^{t-2}\sum_{a=0}^{5k-4t+2}s_{(4^a,3^{k+1-a},2^{4t-4k-1},1^{5k-4t-a+3})},\\[5pt]
D_4(t)=&\sum_{k=0}^{t-2}\sum_{a=0}^{3k-2t}s_{(4^a,3^{2t-k-a-1},2^{3},1^{3k-2t-a+1})},\\[5pt]
D_5(t)=&\sum_{k=0}^{t-1}\sum_{a=0}^{3k-2t-1}s_{(4^a,3^{2t-a-k},2^{2},1^{3k-2t-a})},\\[5pt]
D_6(t)=&\sum_{k=0}^{t-1}\sum_{a=0}^{3k-2t-3}s_{(4^a,3^{2t-a-k},2^{3},1^{3k-2t-a-2})},\\[5pt]
D_7(t)=&\sum_{k=0}^{t-1}\sum_{j=0}^{k-2}\sum_{a=\max(0,k+2i-2t)}^{k+2i-2t}s_{(4^a,3^{2t-k-a},2^{k-i+2})},\\[5pt]
D_8(t)=&\sum_{k=1}^{t-1}\sum_{i=0}^{k-1}\sum_{a=\max(0,k-2t+2i)}^{k-2t+2i}s_{(4^{a},3^{2t-a-k},2^{k-i+2})},\\[5pt]
D_9(t)=&\sum_{k=1}^{t-1}\sum_{j=0}^{2t-2k-2}\sum_{a=\max(0,5k+j-4t+1)}^{5k+j-4t+1}s_{(4^{a},3^{k+1-a+j},2^{4t-4k-2j})}.
\end{align*}

The following lemma shows that one can group the terms of $N_i(t)$
and $M_i(t)$ to reduce the involved triple sums to double sums.

\begin{lem}\label{nm2cd}
For any $t\geq 0$, we have
\begin{align*}
N_2(t)-N_1(t)=&C_1(t)-C_2(t),\\[5pt]
N_4(t)-N_3(t)=&C_3(t),\\[5pt]
N_6(t)-N_5(t)=&C_4(t)+C_5(t),\\[5pt]
M_2(t)-M_1(t)=&D_1(t)-D_2(t)-D_3(t)-D_4(t),\\[5pt]
M_4(t)-M_3(t)=&D_5(t) +D_6(t)-D_7(t),\\[5pt]
M_6(t)-M_5(t)=&D_8(t)-D_9(t).
\end{align*}
\end{lem}

\pf We have
\begin{align*}
N_2(t)-N_1(t)=&
\sum_{k=0}^{t-1}\sum_{j=0}^{2t-2k-2}\sum_{a=0}^{5k+j-4t+2}
s_{(4^a,3^{k-a+j+1},2^{4t-4k-2j-1},1^{5k+j+3-4t-a})}\\[6pt]
&-\sum_{k=0}^{t-1}\sum_{j=0}^{2t-2k-1}\sum_{a=0}^{5k+j-4t-1}
s_{(4^a,3^{k-a+j+1},2^{4t-4k-2j+1},1^{5k+j-1-4t-a})}\\[6pt]
=& \sum_{k=1}^{t}\sum_{j=0}^{2t-2k}\sum_{a=0}^{5k+j-4t-3}
s_{(4^a,3^{k-a+j},2^{4t-4k-2j+3},1^{5k+j-2-4t-a})}\\[6pt]
&-\sum_{k=0}^{t-1}\sum_{j=0}^{2t-2k-1}\sum_{a=0}^{5k+j-4t-1}
s_{(4^a,3^{k-a+j+1},2^{4t-4k-2j+1},1^{5k+j-1-4t-a})}\\[6pt]
=& \sum_{k=1}^{t}\sum_{j=-1}^{2t-2k-1}\sum_{a=0}^{5k+j-4t-2}
s_{(4^a,3^{k-a+j+1},2^{4t-4k-2j+1},1^{5k+j-1-4t-a})}\\[6pt]
&-\sum_{k=0}^{t-1}\sum_{j=0}^{2t-2k-1}\sum_{a=0}^{5k+j-4t-1}
s_{(4^a,3^{k-a+j+1},2^{4t-4k-2j+1},1^{5k+j-1-4t-a})}\\[6pt]
=& \sum_{k=1}^{t}\sum_{j=-1}^{-1}\sum_{a=0}^{5k+j-4t-2}
s_{(4^a,3^{k-a+j+1},2^{4t-4k-2j+1},1^{5k+j-1-4t-a})}\\[6pt]
&+ \sum_{k=1}^{t}\sum_{j=0}^{2t-2k-1}\sum_{a=0}^{5k+j-4t-2}
s_{(4^a,3^{k-a+j+1},2^{4t-4k-2j+1},1^{5k+j-1-4t-a})}\\[6pt]
&-\sum_{k=0}^{t-1}\sum_{j=0}^{2t-2k-1}\sum_{a=0}^{5k+j-4t-1}
s_{(4^a,3^{k-a+j+1},2^{4t-4k-2j+1},1^{5k+j-1-4t-a})}.
\end{align*}
It follows that
$$
N_2(t)-N_1(t)=C_1(t)-C_2(t).
$$
Similarly, one can prove the other identities. This completes the
proof. \qed

It is still necessary to compute $T_9(t+1)-\Delta^{(2,2)}T_9(t)$ and
$T_{10}(t+1)-\Delta^{(2,2)}T_{10}(t)$. Let
\begin{align*}
E_1(t) & = \sum_{a=0}^{t+1}s_{(4^a,3^{t+1-a},1^{t+1-a})},\\[5pt]
E_2(t) & = \sum_{k=0}^t\sum_{a=0}^{3k-2t-1}s_{(4^a,3^{2t-k-a+1},2,1^{3k-2t-1-a})},\\[5pt]
E_3(t) & =
\sum_{k=0}^t\sum_{a=0}^{3k-2t-2}s_{(4^a,3^{2t-k-a+2},1^{3k-2t-2-a})},
\\[5pt]
E_4(t) & = \sum_{a=0}^t
s_{(4^a,3^{t+1-a},1^{t+1-a})},\\[5pt]
E_5(t) & =\sum_{a=0}^{t-1} s_{(4^a,3^{t+1-a},2,1^{t-1-a})},\\[5pt]
E_6(t) & = \sum_{k=0}^{t-1}
\sum_{a=0}^{3k-2t+1}s_{(4^a,3^{2t-k-a},2,1^{3k+2-2t-a})},\\[5pt]
E_7(t) & = \sum_{k=0}^{t-1}
\sum_{a=0}^{3k-2t}s_{(4^a,3^{2t+1-k-a},1^{3k+1-2t-a})},\\[5pt]
E_8(t) & = \sum_{k=0}^{t-1}
\sum_{a=0}^{3k-2t}s_{(4^a,3^{2t-k-a},2,2,1^{3k-2t-a})}\\[5pt]
E_9(t) & = \sum_{k=0}^{t-1}
\sum_{a=0}^{3k-2t-1}s_{(4^a,3^{2t+1-k-a},2,1^{3k-1-2t-a})}.
\end{align*}

The following lemma gives the Schur expansions of
$T_9(t+1)-\Delta^{(2,2)}T_9(t)$ and
$T_{10}(t+1)-\Delta^{(2,2)}T_{10}(t)$.

\begin{lem}\label{t9-t10}
For any $t\geq 0$, we have
\begin{align}
T_9(t+1)-\Delta^{(2,2)}T_9(t)&=E_1(t)+E_2(t)+E_3(t), \label{t9}\\[5pt]
T_{10}(t+1)-\Delta^{(2,2)}T_{10}(t)&=E_4(t)+E_5(t)+E_6(t)+E_7(t)+E_8(t)+E_9(t).\label{t10}
\end{align}
\end{lem}

\pf We will present the proof of the identity (\ref{t10}), because
it is easier to prove (\ref{t9}) also by using Pieri's rule. Recall
that
$$T_{10}(t)= \sum_{k=0}^{t-1}e_ks_{(3^{k+1},2^{2t-2k-2},1)}.$$ We
have
\begin{align*}
T_{10}(t+1) - \Delta^{(2,2)}T_{10}(t)=&
\sum_{k=0}^{t}e_ks_{(3^{k+1},2^{2t-2k},1)}-\Delta^{(2,2)}\left(\sum_{k=0}^{t-1}e_ks_{(3^{k+1},2^{2t-2k-2},1)}\right)\\[6pt]
=& e_ts_{(3^{k+1},1)}+\sum_{k=0}^{t-1}e_ks_{(3^{k+1},2^{2t-2k},1)}\\[6pt]
 & -\sum_{k=0}^{t-1}\Delta^{(2,2)}(e_ks_{(3^{k+1},2^{2t-2k-2},1)}).
\end{align*}
Using the dual version of Pieri's rule, we obtain that
$$e_ts_{(3^{k+1},1)}=\sum_{a=0}^t
s_{(4^a,3^{t+1-a},1^{t+1-a})}+\sum_{a=0}^{t-1}
s_{(4^a,3^{t+1-a},2,1^{t-1-a})},$$ and, for $0\leq k\leq t-1$,
$$e_ks_{(3^{k+1},2^{2t-2k-2},1)}=\sum_{a,b,c,d}s_{(4^a,3^{k+1-a+b},2^{2t-2k-2-b+c},1^{1-c+d})},$$
where
$$0\leq a\leq k+1,\, 0\leq b\leq 2t-2k-2,\, 0\leq c\leq 1,\, d\geq 0,\, a+b+c+d=k.$$
Therefore,
\begin{align*}
e_ks_{(3^{k+1},2^{2t-2k},1)}-\Delta^{(2,2)}(e_ks_{(3^{k+1},2^{2t-2k-2},1)})=\sum_{a,b,c,d}s_{(4^a,3^{k+1-a+b},2^{2t-2k-b+c},1^{1-c+d})},
\end{align*}
where
$$0\leq a\leq k+1,\, 2t-2k-1\leq b\leq 2t-2k,\, 0\leq c\leq 1,\, d\geq 0,\, a+b+c+d=k.$$
In view of the ranges of $b$ and $c$, the above sum is divided into
four parts:
$$
\begin{array}{c}
\sum_{a=0}^{3k-2t+1}s_{(4^a,3^{2t-k-a},2,1^{3k+2-2t-a})}+
\sum_{a=0}^{3k-2t}s_{(4^a,3^{2t+1-k-a},1^{3k+1-2t-a})}\\[8pt]
\rule{16pt}{0pt}+\sum_{a=0}^{3k-2t}s_{(4^a,3^{2t-k-a},2,2,1^{3k-2t-a})}+
\sum_{a=0}^{3k-2t-1}s_{(4^a,3^{2t+1-k-a},2,1^{3k-1-2t-a})}.
\end{array}
$$
This completes the proof. \qed

We are now in a position to give a recurrence relation for $L(r)$ by
using the above lemmas. Note that it is easy to establish a
recurrence relation for $R(r)$. For $t\geq 0$, let
$$
\begin{array}{rcl}
R_{o,1}(t) & = & R(2t+1)-\Delta^{(1,1)}R(2t),\\[5pt]
R_{e,1}(t) & = & R(2t+2)-\Delta^{(1,1)}R(2t+1),\\[5pt]
L_{o,1}(t) & = & L(2t+1)-\Delta^{(1,1)}L(2t),\\[5pt]
L_{e,1}(t) & = & L(2t+2)-\Delta^{(1,1)}L(2t+1),
\end{array}
$$
and
$$
\begin{array}{rcl}
R_{o,2}(t) & = & R_{o,1}(t+1)-\Delta^{(2,2)}R_{o,1}(t),\\[5pt]
R_{e,2}(t) & = & R_{e,1}(t+1)-\Delta^{(2,2)}R_{e,1}(t),\\[5pt]
L_{o,2}(t) & = & L_{o,1}(t+1)-\Delta^{(2,2)}L_{o,1}(t),\\[5pt]
L_{e,2}(t) & = & L_{e,1}(t+1)-\Delta^{(2,2)}L_{e,1}(t).
\end{array}
$$
The following lemma gives the recurrence relations of $L(r)$ and
$R(r)$.

\begin{lem}\label{2-lev}
Let $\mathrm{Par}_{\{3,4\}}(n)$ denote the set of partitions of $n$
whose parts belong to $\{3,4\}$. Then for any $t\geq 0$ we have
\begin{align}
R_{o,2}(t)=\sum_{\lambda\in \mathrm{Par}_{\{3,4\}}(4r+4)}s_\lambda,
\quad
R_{e,2}(t)=\sum_{\lambda\in \mathrm{Par}_{\{3,4\}}(4r+6)}s_\lambda,\label{right2}\\
L_{o,2}(t)=\sum_{\lambda\in \mathrm{Par}_{\{3,4\}}(4r+4)}s_\lambda,
\quad L_{e,2}(t)=\sum_{\lambda\in
\mathrm{Par}_{\{3,4\}}(4r+6)}s_\lambda.\label{left2}
\end{align}
\end{lem}
\pf It is easy to verify the identities in \eqref{right2}. It
remains to prove the identities in \eqref{left2}. Here we will
consider only the identity concerning $L_{o,2}(t)$, since the
identity for $L_{e,2}(t)$ can be justified by the same argument.

Note that by Lemma \ref{lem-ten}, we obtain
\begin{align*} L_{o,2}(t) =&
\sum_{m=1}^4(T_m(t+1)-\Delta^{(2,2)}T_m(t))-
\sum_{m=5}^8(T_m(t+1)-\Delta^{(2,2)}T_m(t))\\[5pt]
&
+(T_9(t+1)-\Delta^{(2,2)}T_9(t))-(T_{10}(t+1)-\Delta^{(2,2)}T_{10}(t)).
\end{align*}

By Lemma \ref{T1-T8}, we find
\begin{align*}
L_{o,2}(t)
=&(T_{41}(t)-T_{61}(t))+(T_{31}(t)-T_{71}(t))+(T_{11}(t)-T_{81}(t))\\[5pt]
&+(T_{21}(t)-T_{51}(t))
+(T_{32}-T_{72}(t))+(T_{42}(t)-T_{62}(t))\\[5pt]
&+(T_{12}(t)-T_{82}(t))+(T_{22}(t)-T_{52}(t))\\[5pt]
&
+(T_9(t+1)-\Delta^{(2,2)}T_9(t))-(T_{10}(t+1)-\Delta^{(2,2)}T_{10}(t)).
\end{align*}
Applying Lemmas \ref{t2nm}-\ref{t9-t10}, we obtain
\begin{align*}
L_{o,2}(t) =&(N_2(t)-N_1(t))+(N_4(t)-N_3(t))+(N_6(t)-N_5(t))\\[5pt]
&+(M_2(t)-M_1(t))+(M_4(t)-M_3(t))+(M_6(t)-M_5(t))\\[5pt]
&
+(E_1(t)+E_2(t)+E_3(t))\\[5pt]
&-(E_4(t)+E_5(t)+E_6(t)+E_7(t)+E_8(t)+E_9(t))\\[5pt]
=&(C_1(t)-C_2(t)+C_3(t)+C_4(t)+C_5(t))\\[5pt]
&+(D_1(t)-D_2(t)-D_3(t)-D_4(t)+D_5(t)\\[5pt]
&+D_6(t)-D_7(t)+D_8(t)-D_9(t))\\[5pt]
& +(E_1(t)+E_2(t)+E_3(t)-E_4(t)-E_5(t)\\[5pt]
&-E_6(t)-E_7(t)-E_8(t)-E_9(t))\\[5pt]
=&[(C_3(t)-E_9(t))+(E_2(t)-E_6(t))-E_5(t)]\\[5pt]
&+[(C_5(t)-D_9(t))+(D_5(t)-E_8(t))+C_4(t)]\\[5pt]
&+[(D_6(t)-D_4(t))+(C_1(t)-D_3(t))-D_2(t)]\\[5pt]
&+[(D_8(t)-D_7(t))+(D_1(t)-C_2(t))]\\[5pt]
&+[(E_1(t)-E_4(t))+(E_3(t)-E_7(t))]\\[5pt]
=&0+0+0+0+\sum_{k=1}^{t+1}s_{(4^{3k-2t-2},3^{4t-4k+4})},
\end{align*}
where the last equality comes from the following relations:
\begin{align*}
C_3(t)-E_9(t)=&\sum_{a=0}^{t-2}s_{(4^a,3^{t-a+1},2,1^{t-a-1})}\\
&-\sum_{k=0}^{t-2}\sum_{a=\max(0,3k-2t+2)}^{3k-2t+2}s_{(4^a,3^{2t-k-a},2)},\\[5pt]
C_5(t)-D_9(t)=&\sum_{k=1}^{t-1}\sum_{a=\max(0,3k-2t)}^{3k-2t}s_{(4^a,3^{2t-k-a},2^2)}\\
&-\sum_{k=1}^{t-1}\sum_{a=\max(0,5k-4t+1)}^{5k-4t+1}s_{(4^a,3^{k+1-a},2^{4t-4k})},\\
C_1(t)-D_3(t)=&\sum_{a=0}^{t-3}s_{(4^a,3^{t-a},2^3,1^{t-a-2})},\\
D_5(t)-E_8(t)=&-\sum_{k=0}^{t-1}\sum_{a=\max(0,3k-2t)}^{3k-2t}s_{(4^a,3^{2t-k-a},2^2)},\\[5pt]
D_6(t)-D_4(t)=&0,\\[5pt]
D_8(t)-D_7(t)=&\sum_{k=1}^{t-1}\sum_{a=\max(0,3k-2t-2)}^{3k-2t-2}s_{(4^a,3^{2t-a-k},2^3)},\\[5pt]
D_1(t)-C_2(t)=&-\sum_{k=0}^{t-2}\sum_{a=\max(0,3k-2t+1)}^{3k-2t+1}s_{(4^a,3^{2t-a-k-1},2^3)},\\[5pt]
E_1(t)-E_4(t)=&s_{(4^{t+1})},\\
E_2(t)-E_6(t)=&\sum_{k=1}^{t}\sum_{a=\max(0,3k-2t-1)}^{3k-2t-1}s_{(4^a,3^{2t-k-a+1},2)},\\[5pt]
E_{3}(t)-E_{7}(t)=&\sum_{k=1}^{t}\sum_{a=\max(0,3k-2t-2)}^{3k-2t-2}s_{(4^a,3^{2t-k-a+2})}.
\end{align*}
So we have obtained the desired Schur expansion of $L_{o,2}(t)$.
This completes the proof. \qed

Based on the above lemma, we obtain the following relations.

\begin{lem}\label{1-lev}
For $t\geq 0$, we have
\begin{align*}
L_{o,1}(t)=R_{o,1}(t),\quad L_{e,1}(t)=R_{e,1}(t).
\end{align*}
\end{lem}
\pf We conduct induction on $t$. We first consider that relation
$L_{o,1}(t)=R_{o,1}(t)$. Clearly, the equality holds for $t=0$ and
$t=1$. Assume that $L_{o,1}(t-1)=R_{o,1}(t-1)$. Note that
\begin{align*}
L_{o,1}(t)=L_{o,2}(t-1)+\Delta^{(2,2)}L_{o,1}(t-1),\\[3pt]
R_{o,1}(t)=R_{o,2}(t-1)+\Delta^{(2,2)}R_{o,1}(t-1).
\end{align*}
It follows  from Lemma \ref{2-lev} that $L_{o,2}(t-1)=R_{o,2}(t-1)$.
Hence by induction we have $L_{o,1}(t)=R_{o,1}(t)$. Similarly, it
can be shown that $L_{e,1}(t)=R_{e,1}(t)$. \qed

 We have established
the recurrence relations  for $L(r)$ and $R(r)$. We now proceed to
prove Theorem \ref{condc1} based on these recurrence relations.

\noindent \textit{Proof of Theorem \ref{condc1}.}  We conduct
induction on $r$. It is easy to verify that the identity holds for
$r=1$. Assume that $L(r)=R(r)$. We proceed to prove that
$L(r+1)=R(r+1)$. If $r=2t$, then
\begin{align*}
L(r+1)=L(2t+1)&=\Delta^{(1,1)}L(2t)+L_{o,1}(t),\\
R(r+1)=R(2t+1)&=\Delta^{(1,1)}R(2t)+R_{o,1}(t).
\end{align*}
From Lemma \ref{1-lev} it follows that $L(r+1)=R(r+1)$. If $r=2t+1$,
then
\begin{align*}
L(r+1)=L(2t+2)&=\Delta^{(1,1)}L(2t+1)+L_{e,1}(t),\\
R(r+1)=R(2t+2)&=\Delta^{(1,1)}R(2t+1)+R_{e,1}(t).
\end{align*}
We also have $L(r+1)=R(r+1)$ by Lemma \ref{1-lev}. This completes
the proof. \qed

\section{$q$-log-convexity}

In this section, we aim to prove Theorem \ref{Liuwang1} and Theorem
\ref{Liuwang2}. The proof of Theorem \ref{Liuwang1} is based on the
 identity \eqref{main-eq} in the previous section.

\noindent \textit{Proof of Theorem \ref{Liuwang1}.} For any $n\geq
1$, we have $\ps_n^1(e_k)={n\choose k}$. So we see that
$$W_n(q)=\sum_{k=0}^n{n\choose k}^2q^k=\sum_{k=0}^n(\ps_n^1(e_k))^2q^k.$$
Thus, for any $r\geq 0$, the coefficient of $q^{r}$ in
$W_{n-1}(q)W_{n+1}(q)-(W_n(q))^2$ is given by
$$\sum_{k=0}^r{\ps_{n-1}^1(e_k)}^2{\ps_{n+1}^1(e_{r-k})}^2-{\ps_{n}^1(e_k)}^2
{\ps_{n}^1(e_{r-k})}^2.$$ By \eqref{decOrder}, the above sum equals
$$\ps_{n-1}^1\left(\sum_{k=0}^r{e_k}^2{(e_{r-k}+2e_{r-k-1}+e_{r-k-2})}^2
-{(e_k+e_{k-1})}^2 (e_{r-k}+e_{r-k-1})^2\right).$$ To evaluate the
above sum, we first expand  the squares, and then apply the
following relations
\begin{align*}
&\sum_{k=0}^r e_k^2e_{r-k-2}^2=\sum_{k=0}^r
e_{k-1}^2e_{r-k-1}^2,\\[3pt]
&\sum_{k=0}^r e_{k-1}^2e_{r-k}^2=\sum_{k=0}^r e_{k}^2e_{r-k-1}^2,\\[3pt]
&\sum_{k=0}^r e_{k}^2e_{r-k}e_{r-k-1}= \sum_{k=0}^r
e_{r-k}^2e_{k}e_{k-1},\\[3pt]
&\sum_{k=0}^r e_{k}^2e_{r-k}e_{r-k-2}=
\sum_{k=0}^re_{r-k}^2e_{k}e_{k-2},\\[3pt]
&\sum_{k=0}^r e_{k-1}^2e_{r-k}e_{r-k-1}=\sum_{k=0}^r
e_{r-k-1}^2e_{k}e_{k-1} =\sum_{k=0}^r e_{k}^2e_{r-k-1}e_{r-k-2}.
\end{align*}
After simplification we arrive at the following expression
\begin{align*}
&2\,\ps_{n-1}^1\left(\sum_{k=0}^re_{k-1}^2e_{r-k}^2+e_{k-2}e_{k}e_{r-k}^2-2e_{k-1}e_{k}e_{r-k-1}e_{r-k}\right).
\end{align*}
By \eqref{main-eq}, we see that the sum in the above expression is
Schur positive.  This implies that the polynomials $W_n(q)$ form a
$q$-log-convex sequence, completing the proof. \qed

To prove Theorem \ref{Liuwang2}, we first introduce the following
polynomials. For any $n\geq 1$ and $0\leq r\leq 2n$, define
\begin{eqnarray*}
f_1(x) & = & (n+1)^2(n-x+1)^2(n-x)^2,\\[5pt]
f_2(x) & = & (n+1)^2(n-(r-x)+1)^2(n-(r-x))^2,\\[5pt]
f_3(x)  & = & n^2(n-x+1)^2(n-(r-x)+1)^2.
\end{eqnarray*}
Let $$f(x)=f_1(x)+f_2(x)-2f_3(x).$$

\noindent \textit{Proof of Theorem \ref{Liuwang2}.} It suffices to
prove that the polynomials $W_n(q)$ satisfy the  conditions in
Theorem \ref{Liuwang3}. Clearly, for any $n$ and $r$, if $k\leq
r-n-1$, then $n\leq (r-k)-1$ and $\alpha(n,r,k)=0$. We only need to
determine the sign of $\alpha(n,r,k)$ for $r-n-1<k\leq \lfloor
\frac{r}{2}\rfloor$. It is easy to see that $\alpha(n,r,k)$ can be
rewritten as
\begin{equation}\label{tranform}
\alpha(n,r,k)=\frac{1}{n^2(n-k+1)^2(n-r+k+1)^2}{n\choose
k}^2{n\choose {r-k}}^2f(k).
\end{equation}
Let us consider the value of $f(k)$ for a given $r$. Taking the
derivative of $f(x)$ with respect to $x$, we obtain that
 $$f'(x)=2(2x-r)g(x),$$
 where
\begin{align*}
g(x)&=2x^2+4nx^2-4nrx-2rx+4nr^2-17n^2r+2n^2r^2-12nr\\[5pt]
&\rule{12pt}{0pt}-8n^3r+1+8n+21n^2+8n^4+22n^3+2r^2-3r.
\end{align*}
Differentiating with respect to $x$, we find that
\begin{equation*}
g'(x)=2(2x-r)(1+2n).
\end{equation*}
This yields that $g'(x)\leq 0$ for $x\leq \frac{r}{2}$. Thus $g(x)$
is decreasing on the interval $(-\infty,\frac{r}{2}]$. Since
$g(x)\rightarrow +\infty$ when $x\rightarrow -\infty$, the
polynomial $g(x)$ has at most one real root, say $x_0$ if it exists,
on the interval $(-\infty,\frac{r}{2}]$. Hence, either $f'(x)$ has a
unique root $\frac{r}{2}$ on the interval $(-\infty,\frac{r}{2}]$,
or it has two real roots $x_0$ and $\frac{r}{2}$. In the former
case, $f(x)$ is decreasing on the interval $(-\infty,\frac{r}{2}]$.
In the latter case, $f(x)$ is decreasing on the interval
$(-\infty,x_0]$ and increasing on the interval $[x_0,\frac{r}{2}]$.
Combining the two cases, it suffices to show that $f(r/2)\leq 0$.
Observe that
\begin{align*}
f(r/2) & =
2(n+1)^2(n-\frac{r}{2}+1)^2(n-\frac{r}{2})^2-2n^2(n-\frac{r}{2}+1)^4\\[6pt]
  & =  -\frac{r(2n(2n-r)+(2n-r)+2n)(2+2n-r)^2}{8},
\end{align*}
which is nonpositive for $0\leq r\leq 2n$. By \eqref{tranform},
there exists an integer $k'=k'(n,r)$ such that $\alpha(n,r,k)\geq 0$
for $k\leq k'$ and $\alpha(n,r,k)\leq 0$ for $k> k'$. Therefore, by
Theorem \ref{Liuwang1} and the conditions of Liu and Wang, we deduce
that the linear transformation defined by the triangular array $\{{n
\choose k}^2\}_{0\leq k\leq n}$  preserves  $q$-log-convexity. This
completes the proof. \qed

\vskip 8pt

\noindent {\bf Acknowledgments.} This work was supported by  the 973
Project, the PCSIRT Project of the Ministry of Education, the
Ministry of Science and Technology, and the National Science
Foundation of China.

\end{document}